\documentclass[12pt]{article}
\usepackage{amssymb}
\usepackage{hyperref}

\def\al{\alpha}
\def\be{\beta}
\def\cat(#1){{\hat{\phantom{a}}\langle{#1}\rangle}}
\def\df{\em}
\def\diagonal{\Delta}

\def\la{\langle}
\def\om{\omega}
\def\proof{\par\noindent Proof\par\noindent}
\def\pr{\prime}
\def\qed{\par\noindent QED\par\bigskip}
\def\rank{{\rm rank}}
\def\ra{\rangle}
\def\res{\upharpoonright}
\def\rmand{\mbox{ and }}
\def\rmiff{\mbox{ iff }}
\def\rmor{\mbox{ or }}

\def\si{\sigma}
\def\sm{{\setminus}}
\def\su{\subseteq}
\def\tless{{\,\lhd\,}}

\newtheorem{theorem}{Theorem}
\newtheorem{define}[theorem]{Definition}
\newtheorem{lemma}[theorem]{Lemma}
\newtheorem{prop}[theorem]{Proposition}

\def\address{\begin{flushleft}
Arnold W. Miller \\
miller@math.wisc.edu \\
http://www.math.wisc.edu/$\sim$miller\\
University of Wisconsin-Madison \\
Department of Mathematics, Van Vleck Hall \\
480 Lincoln Drive \\
Madison, Wisconsin 53706-1388 \\
\end{flushleft}}

\begin{document}

\begin{center}
{\large A hierarchy of clopen graphs on the Baire space}
\end{center}

\begin{flushright}
Arnold W. Miller\\
Oct 2012
\end{flushright}

We say that $E\su X\times X$ 
is a {\df clopen graph} on $X$ iff 
$E$ is symmetric and irreflexive and clopen relative to
$X^2\sm \diagonal$ where $\diagonal=\{(x,x):x\in X\}$ is the diagonal. 
Equivalently $E\su [X]^2$ and for all $x\neq y\in X$ there
are open neighborhoods  $x\in U$ and $y\in V$ such that
either $U\times V\su E$ or $U\times V\su X^2\sm E$.

For clopen graphs $E_1,E_2$ on spaces $X_1,X_2$, we say that
$E_1$ {\df continuously reduces} to $E_2$ iff there is a continuous
map $f:X_1\to X_2$ such that for every
$x,y\in X_1$
$$(x,y)\in E_1\rmiff (f(x),f(y))\in E_2.$$
Note that $f$ need not be one-to-one but there should be no
edges in the preimage of a point.
If $f$ is a homeomorphism to its image,
then we say that $E_1$ {\df continuously embeds} into $E_2$.

\begin{theorem}
There does not exist countably 
many clopen graphs on the Baire space, $\om^\om$, such that every clopen graph 
on $\om^\om$ can be continuously
reduced to one of them.  However, there are $\om_1$ clopen graphs 
on $\om^\om$
such that every clopen graph on $\om^\om$
continuously embeds into one of them.
\end{theorem}

Since one can take a countable clopen separated union of
countably many clopen graphs on $\om^\om$, having countably many
is the same as having one universal graph.

\begin{define}
For $R\su \om^\om\times\om^\om$, $C$ and $D$ clopen
subsets of $\om^\om$, and $\al$ an ordinal define
\begin{enumerate}
 \item $\rank_R(C\times D)=0$ iff $C\times D\su R$ or
$R\cap (C\times D)=\emptyset$
 \item $\rank_R(C\times D)\leq \al$ iff there are partitions
of $C$ and $D$ into clopen sets:
 $C=\sqcup_{i<\om}C_i$ and $D=\sqcup_{j<\om}D_j$
such that $\rank_R(C_i\times D_j)<\al$ for all $i,j$ in $\om$.
\end{enumerate}
We use $\sqcup$ to mean disjoint union.
\end{define}

\bigskip
\noindent Since we allow $C_i$'s and $D_j$'s to be empty, it is clear that:

\begin{prop}\label{propclear}
If $\rank_R(C\times D)\leq \al$ and $C^\pr\su C$ and $D^\pr\su D$,
then $\rank_R(C^\pr\times D^\pr)\leq \al$.
\end{prop}

\noindent More generally:

\begin{prop}\label{propreduction}
Suppose $f:C \sqcup D \to C^\pr\sqcup D^\pr$ is a continuous
reduction of $R \su C \times D $ to $R^\pr \su C^\pr \times D^\pr$
and $f^{-1}(C^\pr)=C $
and $f^{-1}(D^\pr)=D $.  Then
$$\rank_{R }(C \times D )\leq \rank_{R^\pr}(C^\pr\times D^\pr).$$
\end{prop}
\proof
Since $f$ is continuous, clopen partitions 
$C^\pr=\sqcup C^\pr_i$ and $D^\pr=\sqcup D^\pr_j$ 
induce clopen partitions
$C=\sqcup_if^{-1}(C^\pr_i)$ and $D=\sqcup_jf^{-1}(D^\pr_j)$.
\qed

\begin{define}
For $E$ a clopen graph on $\om^\om$ define
$$\rank(E)=\sup\{\rank_E(C\times D)\;:\; 
\mbox{ $C$ and $D$ are disjoint clopen sets}\}.$$
\end{define}

\begin{lemma} \label{four}
If $E$ is a clopen graph on $\om^\om$,
then $\rank(E)<\om_1$.
\end{lemma}
\proof
Given incomparable $s_0,t_0\in \om^{<\om}$ with
the same length look at the tree $T$:

\bigskip\noindent $(s,t)\in T$ iff
\begin{enumerate}
\item $s_0\su s$, $t_0\su t$, $|s|=|t|$, and
\item both $([s]\times [t])\;\cap E$ and 
$([s]\times [t])\;\sm E$ are nonempty.
\end{enumerate}

Let
$$T^*=\{(s_0,t_0)\}\cup\{(s\cat(i),t\cat(j))\;:\; (s,t)\in T \rmand i,j\in\om\}.$$ 
Since $E\cap([s_0]\times [t_0])$ is clopen, $T^*$ is well-founded and $T^*\sm T$ is the
set of the terminal nodes of $T^*$.  Let $r$ be
the standard
rank function on $T^*$, i.e., 
\begin{itemize}
\item $r(s,t)=0$ iff $(s,t)\in T^*\sm T$ 
\item $r(s,t)=\sup\{r(s\cat(i),t\cat(j))+1\;:\;i,j\in\om\}\;\;\;\;$
if $(s,t)\in T$.
\end{itemize}
Note that $\rank_E([s]\times [t])\leq r(s,t)$ for $(s,t)\in T^*$.
Take any countable ordinal $\al$ such that for
every $s\in \om^{<\om}$ and
distinct $i,j\in\om$ we have that $\rank_E([s\cat(i)]\times [s\cat(j)])\leq \al$.

For any $s,t\in\om^{<\om}$ which are incomparable, let $n$
be the least such that $s\res n\neq t\res n$.  Then $[s]\su[s\res n]$
and $[t]\su[t\res n]$ and so by Proposition \ref{propclear},
$\rank_E([s]\times [t])\leq \al$.
But any nonempty open set $U$ can be written as pairwise
disjoint basic clopen sets, i.e., $U=\sqcup_{i<\om} [s_i]$
where $s_i\in \om^{<\om}$.
(To see this just take for any $x\in U$ the least $n$ with
$[x\res n]\su U$.)
Hence for any disjoint clopen sets $C,D$ we have that
$\rank_E(C\times D)\leq \al+1$.
And so $\rank(E)\leq\al+1$.
\qed

\begin{lemma}\label{no-universal}
For any $\al<\om_1$ there exists a clopen graph $E_\al$ on $\om^\om$
such that if $E$ is any
clopen graph on the Baire space such that $E_\al$ 
is continuously reducible to $E$,
then $\rank(E)\geq \al$.
\end{lemma}
\proof
Let $Q=\{in,out\}$ and $\al$ any countable limit ordinal.  Put
$\Gamma_\al=\om\times(Q\cup \al)$ with the discrete topology and
define a clopen relation
$R_\al\su \om^\om\times \Gamma_\al^\om$
as follows.
Given $x\in\om^\om$ and $y\in \Gamma_\al^\om$ construct
sequences $m_i,n_i\in\om$ and $\al_i\in \al\cup Q$ as 
follows.
\begin{itemize}
\item $x(0)=m_0$ and $y(m_0)=(n_0,\al_0)$
\item $x(n_{i-1})=m_{i}$ and $y(m_i)=(n_i,\al_i)$ for $i\geq 1$.
\end{itemize}
To determine whether or not $(x,y)\in R_\al$ look at the first
$i$ such that either $\al_i\in Q$ or ($i>0$ and 
$\al_i\in\al$ but {\bf not} $\al_i<\al_{i-1}$).
Note that such an $i$ must always
occur since otherwise we would get an infinite descending sequence
of ordinals.  Let $i_0$ be the first such $i$ and
put $(x,y)\in R_\al$ iff $\al_{i_0}=in$.

Note that $R_\al$ is clopen since given any
$x,y\in\om^\om$ we can choose $N$ sufficiently large so that
every pair in $[x\res N]\times [y\res N]$ will terminate the same
way $(x,y)$ did.

\bigskip\noindent {\bf Claim \ref{no-universal}.1}.
Suppose $s\in \om^{<\om}$ and $t\in \Gamma_\al^{<\om}$
have the property 
that we can define the sequences
$m_i$ and $(n_i,\al_i)$ for $i<N$
using the same prescription as above:
\begin{enumerate}
 \item $s(0)=m_0$ and $t(m_0)=(n_0,\al_0)$, 
 \item $s(n_{i-1})=m_{i}$ and $t(m_i)=(n_i,\al_i)$ for $1\leq i<N$,
 \item $\al_0>\al_1>\cdots>\al_{N-1}$ are all ordinals, and
 \item $s(n_{N-1})=m_N$,
 \item however $m_N\geq|t|$ so we have not yet determined $\al_N$.
\end{enumerate}
Then $\rank_{R_\al}([s]\times [t])\geq\al_{N-1}$.

\medskip\proof 
Suppose that
$\al_{N-1}$ is the least ordinal for which this could be false
(for any $s,t,N$)
and let $\rank_{R_\al}([s]\times [t])=\be<\al_{N-1}$.
It is easy to check that if $\al_{N-1}>0$ then $\be$ cannot be
zero since we may find extensions $t_{in},t_{out}$ of $t$ with
$t_{in}(m_N)=(\cdot,in)$ and $t_{out}(m_N)=(\cdot,out)$.

Let $[s]=\sqcup_i C_i$ and 
$[t]=\sqcup_j D_j$ be clopen partitions with
$\rank_{R_\al}(C_i\times D_j)<\be$ for all $i,j$.  Extend
$s\su s^\pr$ so that $[s^\pr]\su C_i$ for some $i$.  Extend
$t$ to $t^\pr$ so that $t^\pr(m_N)=(|s^\pr|,\be)=^{def}(n_N,\al_N)$
and $[t^\pr]\su D_j$ for some $j$.
Finally extend $s^\pr$
by putting $s^{\pr\pr}=s^\pr\cat(|t^\pr|)$ so 
$s^{\pr\pr}(\al_N)=m_N=|t^\pr|$.
Now we are in the same situation as before except $\al_N=\be$
is now defined. But
$$\rank_{R_\al}([s^{\pr\pr}]\times [t^{\pr}])\leq\rank_{R_\al}(C_i\times D_j)
<\be=\al_{N}.$$
This violates the minimality of $\al_{N-1}$ and so proves
the Claim.
\qed

Now for any limit ordinal $\al$ and $\be<\al$ let $|s|=2$ and 
$|t|=1$ be defined by
$s(0)=0,s(1)=1$ and $t(0)=(1,\be)$. By the claim
$\rank_{R_\al}([s]\times [t])\geq\be$ and since these
exist for every $\be<\al$,
it follows that $\rank(R_\al)\geq\al$.

Now we adjust $R_\al$ to make its domain and range disjoint.
Identify $\Gamma_\al$ with $\om$ and define 
$S_\al=\{(0\cat(x),1\cat(y))\;:\; (x,y)\in R_\al\}$.
Then $S_\al\su C\times D$ where 
$C=[\la 0\ra]$ and $D=[\la 1\ra]$ are disjoint clopen sets.
Clearly $\rank(S_\al)\geq\al$ as it is a copy of $R_\al$.
Let $C_1$ and $D_1$ be nonempty clopen sets such that
$C\sqcup C_1\sqcup D\sqcup D_1=\om^\om$.
Let
$$P_\al=S_\al\cup (C_1\times D) \cup (C\times D_1)$$
Since $P_\al\cap (C\times D)=S_\al$ we know $\rank(P_\al)\geq\al$.
If we let $A=C\cup C_1$ and $B=D\cup D_1$ then $A$ and $B$ are 
complementary clopen sets with $P_\al\su A\times B$ and for every $x\in A$
there is a $y\in B$ with $(x,y)\in P_\al$ and for every $y\in B$
there is an $x\in A$ with $(x,y)\in P_\al$.  (This property that
everything is connected to something else might have already been
true of $R_\al$ but if not, in this step we have added it.)

Finally we define the clopen graph $E_\al$.  We put

\centerline{$(x,y)\in E_\al$ iff $(x,y)\in P_\al$ or $(y,x)\in P_\al$.}

\noindent Then $E_\al$ is a true clopen graph, i.e., $E_\al$ is
a clopen relation in $(\om^\om)^2$ which is symmetric
and irreflexive.  Also $\om^\om=A\sqcup B$ where every element
$\om^\om$ is connected to something else, 
but neither $A$ nor $B$ contain two elements
which are connected.

\bigskip\noindent {\bf Claim \ref{no-universal}.2}.
Suppose $\al$ a countable limit ordinal
and there is a continuous reduction of $E_\al$ to
a clopen graph $E$.  Then $\rank(E)\geq\al$.

\proof
Let $f:\om^\om\to\om^\om$ be a continuous reduction and suppose
for contraction that $\rank(E)=\be<\al$.  Then in particular
for every $(x,y)\in A\times B$
$$(x,y)\in P_\al\rmiff (f(x),f(y))\in E.$$

Let $A^\pr=f(A)$ and $B^\pr=f(B)$
We show that not only are these sets disjoint
but they have a stronger separation property.

For every $z\in A^\pr\cup B^\pr$ there exists some $n$ such that 
$f^{-1}([z\res n])\su A$ or $f^{-1}([z\res n])\su B$.
To see why let $z=f(x)$ for some $x\in A$.
By our construction of $P_\al$ there is a $y\in B$
with $(x,y)\in P_\al$.  By the reduction $(f(x),f(y))\in E$
and since $E$ is clopen
$[f(x)\res n]\times [f(y)\res n]\su E$  for some $n$.  
So if $f(u)\in [f(x)\res n]$
then $(f(u),f(y))\in E$ and
so $(u,y)\in P_\al$.  But this implies $u\in A$ since $y\in B$.

Now define $\Sigma\su \om^{<\om}$ by
$$\Sigma=\{s\in\om^{<\om}\;:\; f^{-1}[s]\su A \rmor f^{-1}[s]\su B\}$$ 
and let
$$\Sigma_0=\{s\in\Sigma\;:\;\forall t\in \Sigma \;\;t\su s \rightarrow
t=s\}.$$
Note that the elements of $\Sigma_0$ are pairwise incomparable
and that 
$$A=\sqcup_{s\in\Sigma_0}\{f^{-1}[s]\;:\; f^{-1}[s]\su A\}$$ 
and
$$B=\sqcup_{t\in\Sigma_0}\{f^{-1}[t]\;:\; f^{-1}[t]\su B\}$$ are 
clopen partitions of $A$ and $B$.  
Since $\rank(E)\leq\be$ for any distinct $s,t\in \Sigma_0$ we
have that $\rank_E([s]\times[t])\leq \be$.  By Proposition \ref{propreduction}
we get that
$\rank(P_\al)\leq \be+1<\al$, which is a contradiction.
This proves Lemma \ref{no-universal}.
\qed

\bigskip

\begin{lemma}\label{univ}
There exists $U_\al$ for $\al<\om_1$ clopen graphs on $\om^\om$ such
that for every clopen graph $E$ on $\om^\om$ there exists $\al<\om_1$ such
that $E$ continuously embeds into $U_\al$.
\end{lemma}
\proof

For any $s\in\om^{<\om}$ except the trivial sequence $\la\ra$ let
$s^*$ be the parent of $s$, i.e., the unique $s^*\su s$ and $|s^*|=|s|-1$.

Let $\al$ be a countable ordinal, $Q=\{in,out\}$ (or
more generally any countable set).  A pair $(T,l)$
is an {\df $\al$-tree} iff $T$ is a subtree of $\om^{<\om}$
and $l:D\to\al\cup Q$ where 
$D=\{\{s,t\}\in [T]^2\;:\; |s|=|t|\rmand s\neq t\}$
and $l$ satisfies:

if $(s,t)\in D$ and $s^*\neq t^*$ then
 \begin{enumerate}
   \item if $l(s^*,t^*)\in\al$ then $l(s,t)<l(s^*,t^*)$ or $l(s,t)\in Q$
   \item if $l(s^*,t^*)\in Q$ then $l(s,t)=l(s^*,t^*)$.
 \end{enumerate}
Note that $l$ is only defined on pairs with $s\neq t$ of the same
length.
Also if $s^*=t^*$, then $l(s,t)$ can be anything in $\al\cup Q$.
A compact way of stating the above two conditions would be by
taking the binary relation $\tless$ on $\al\cup Q$ defined
by $x\tless y$ iff 
\begin{enumerate}
\item $x,y\in\al$ and $x<y$,
\item $x\in Q$ and $y\in \al$, or
\item $x,y\in Q$ and $x=y$.
\end{enumerate}
Then our condition on $l$ is equivalent to:

if $(s,t)\in D$ and $s^*\neq t^*$ then $l(s,t)\tless l(s^*,t^*)$.

\bigskip

Given any clopen graph $E$ 
we describe the canonical $\al$-tree $(\om^{<\om},l)$
associated with it.  For any distinct $s,t$ of the same
length if $[s]\times [t]\su E$, then put $l(s,t)=in$, if
$([s]\times [t])\cap E=\emptyset$, then put $l(s,t)=out$.

Let $P=\{(s,t)\;:\;l(s,t)\in Q\}$
and note that $P$ is closed downward.
For any $s$ and distinct $i,j\in\om$ the tree
$$T_{s,i,j}=\{(t_1,t_2): s\cat(i)\su t_1,\; s\cat(j)\su t_2,\;
|t_1|=|t_2|,\rmand
(t_1,t_2)\notin P\}$$
is a well-founded tree because $E$ is clopen.  
Let $l\res T_{s,i,j}$ be its rank function.
Picking $\al$ large enough makes $(\om^{<\om},l)$ an $\al$-tree.

Next we construct a universal $\al$-tree $(\om^{<\om},L)$.
It will be very strongly universal in the following sense:
Suppose that $(T,l)$ is any $\al$-tree.  Then there
will exists $\si:\om^{<\om}\to\om^{<\om}$ which is tree embedding,
i.e, 
\begin{enumerate}
\item $\si$ is one-to-one and level preserving, i.e.,  $|\si(s)|=|s|$
\item $\si$ preserves the tree ordering, i.e., 
  $s\su t$ implies $\si(s)\su \si(t)$
\item $\si$ preserves the labeling on edges, i.e.,
$l(s,t)=L(\si(s),\si(t))$ for 
any distinct $s,t$ of the same length.
\end{enumerate}
It easy to see that $\si$ induces a continuous embedding
$f:[T]\to \om^\om$ by $f(x)=\bigcup_{n<\om}\si(x\res n)$ 
which reduces the graph associated to $l$ to the one associated
with $L$.

We construct $L$ to have the following property:

\begin{quote}
For any $n<\om$,
$p\in\om^{n}$,
finite $F\su\om^{n+1}$, and
$f:F\to \al\cup Q$ consistent with $L$,
there will be infinitely many $t\in\om^{n+1}$ with
$t^*=p$ such that $L(t,s)=f(s)$ for all $s\in F$.
By $f$ consistent with $L$ we mean:
for all $s\in F$ if $s^*\neq p$, then $f(s)\tless L(p,s^*)$.
\end{quote}

First let us check that it is possible to construct $L$ with this
property.  Let $(p_n,F_n,f_n)$ list with infinitely many
repetitions all triples $(p,F,f)$ with $p\in\om^{<\om}$,
$F\su T\cap \om^{k+1}$ where $|p|=k$, and
$f:F\to \al\cup Q$ arbitrary.   Construct $(T_n,L_n)$ an $\al$-tree
with $T_n$ finite, $T_n\su T_{n+1}$ and $L_n\su L_{n+1}$ and
if $p_n\in T_n$, $F_n\su T_n$, and 
$f_n$ consistent with $L_n$, then there
exists $t\in T_{n+1}\sm T_n$ with 
$p_n=t^*$ such that $L_{n+1}(s,t)=f_n(s)$ for all $s\in F_n$.
This can be done as follows: choose any $t\notin T_n$ with
$t^*=p_n$. For $s\in F_n$ define $L_{n+1}(s,t)=f_n(s)$.  For
all other $s\in T_n$ with $|s|=|t|$ and $s^*\neq t^*$
put $L_{n+1}(s,t)=q$
for any $q\in Q$ with $q\tless L_n(s^*,t^*)$.

Second let us check that this property is all that is needed for universality.
Write any $\al$-tree as an increasing union of finite subtrees $T_n$
gotten by adding one new
child to some node from $T_n$, i.e.,
$T_{n+1}=T_n\cup \{r_n\}$ where $r_n^*\in T_n$ but $r_n\notin T_n$.
The map $\si$ is constructed by extending $\si\res T_n$ to $T_{n+1}$
by defining $\si$ at $r_n$.  Without all the subscripts one step
looks like this:

Suppose $(T,l)$ is a finite $\al$-tree and
$r\in T$ has no child and let $T_0=T\sm\{r\}$.  Suppose
that $\si:T_0\to\om^{<\om}$ is a tree embedding 
of $(T_0,l\res [T_0]^2)$ into $(\om^{<\om},L)$.
Suppose $|r|=n$, $F=\si(T_0)\cap\om^n$, $p=\si(r^*)$, and
$f:F\to\al\cup Q$ is defined by $f(\si(s))=l(r,s)$.
By our property there are infinitely many $t$
such that we can extend $\si$ to $T$ by defining $\si(r)=t$.
This proves the Lemma.\footnote{This type of argument is familiar
to model theorists who would refer to it as joint embedding, amalgamation,
and universal Fraisse structure.  Set theorists would say
its like Cantor's proof that every countable linear order embeds into
the rationals.}
\qed

Theorem 1 follows immediately from the three Lemmas.

\begin{center}
Remarks
\end{center}
Theorem 1 settles a question of Stefan Geschke \cite{ges}.
It was motivated by his result
that the smallest cardinality of a family of
clopen graphs on the Cantor space, $2^\om$, such that every
such graph can be continuously embedded into some member
of the family is exactly
${\mathfrak d}$, the dominating number.  Geschke also showed that
there is a clopen graph on $\om^\om$ universal for all
clopen graphs on $2^\om$.

The family of $U_\al$ in Lemma \ref{univ} are also universal
for all clopen graphs on closed subsets of $\om^\om$ and
hence for all clopen graphs on zero dimensional Polish spaces.

Recall that a clopen graph $E$ on $X$ is {\df true clopen} 
iff $E\su X^2$ is symmetric irreflexive and clopen in
$X^2$ - not just clopen in $X^2\sm \diagonal$.  The proof
of Lemma \ref{no-universal} shows that in fact there is no
clopen graph on $\om^\om$ which is universal for all true
clopen graphs on $\om^\om$.
Note that if $E_1$ is
continuously reducible to $E_2$ and $E_2$ is true
clopen, then $E_1$ is true clopen.  Also if $E$ is
true clopen, then there exists a clopen partition
$\om^\om=\sqcup_{i<\om} C_i$ such that $C_i^2\cap E=\emptyset$
for each $i<\om$.  Using this we can
vary the proof of Lemma \ref{univ} to produce true
clopen $U_\al^\pr$ for $\al<\om_1$ such that every true 
clopen graph continuously
embeds into one of them.  Construct a $\al$-universal tree 
$L^\prime$ similar to $L$ but
satisfying: if $s,t$ are distinct, $|s|=|t|=n>1$,
and $s(0)=t(0)$, then $L^\pr(s,t)=out$.  Hence we are
thinking of replacing $C_i$ with $[\la i\ra]$.

In the case of unary predicates continuous reducibility is called
Wadge reducibility, i.e., for $A,B\su\om^\om$ define
$A\leq_WB$ iff there exists a continuous $f:\om^\om\to\om^\om$
such that $x\in A$ iff $f(x)\in B$.
For a generalization of Wadge reducibility to Borel labellings
in a better-quasi-order
see van Engelen, Miller, and Steel \cite{van}.
Louveau and Saint-Raymond \cite{louv} contains
some results about the quasi-order of Borel linear
orders under embeddability.
Even for finite graphs the $n$-cycles are pairwise incomparable 
under graph embedding, so we don't get a well-quasi-order.
However there are weaker
notions of reducibility under which
finite graphs are well-quasi-ordered,
see Robertson and Dale \cite{robertson}.  Perhaps there
is a natural notion of reducibility for clopen graphs that
gives a well-quasi-ordering.

\bigskip
\address
\end{document}